\documentclass{amsart}

\usepackage{amssymb}
\usepackage[all]{xy}

\newtheorem{thm}{Theorem}

\newtheorem{lem}[thm]{Lemma}
\newtheorem{cor}[thm]{Corollary}

\newtheorem{prop}[thm]{Proposition}

\newtheorem{conj}[thm]{Conjecture}
   
\theoremstyle{definition}
\newtheorem{defn}[thm]{Definition}

\newtheorem{say}[thm]{}
\newtheorem{exmp}[thm]{Example}

\newtheorem{ques}[thm]{Question}    

\newtheorem{rem}[thm]{Remark}          

\newtheorem*{ack}{Acknowledgments}      

\newtheorem{defn-thm}[thm]{Definition--Theorem}  
\newtheorem{defn-lem}[thm]{Definition--Lemma}  

\theoremstyle{remark}

\let \cedilla =\c
\renewcommand{\c}[0]{{\mathbb C}}

\newcommand{\z}[0]{{\mathbb Z}}

\renewcommand{\r}[0]{{\mathbb R}} 

\renewcommand{\a}[0]{{\mathbb A}}

\newcommand{\s}[0]{{\mathbb S}}

\newcommand{\p}[0]{{\mathbb P}}
\newcommand{\f}[0]{{\mathbb F}}
\newcommand{\q}[0]{{\mathbb Q}}
\newcommand{\map}[0]{\dasharrow}
\newcommand{\qtq}[1]{\quad\mbox{#1}\quad}

\newcommand{\aut}[0]{\operatorname{Aut}}

\newcommand{\tsum}[0]{\textstyle{\sum}}

\newcommand{\bir}[0]{\operatorname{Bir}}

\def\into{\DOTSB\lhook\joinrel\to}

\def\loccoh#1.#2.#3.#4.{H^{#1}_{#2}(#3,#4)}

\DeclareMathAlphabet{\mathchanc}{OT1}{pzc}%
                                {m}{it}

\newcommand{\GL}{\mathrm{GL}}
\newcommand{\PGL}{\mathrm{PGL}}

\usepackage[all]{xy}\xyoption{dvips}

\newcommand{\tprod}[0]{\textstyle{\prod}}

\begin{document}
\bibliographystyle{amsalpha}

\today

\title{Algebraic  hypersurfaces}
  \author{J\'anos Koll\'ar}

\begin{abstract} We give an  introduction to the study of algebraic hypersurfaces, focusing on the problem of when two hypersurfaces are isomorphic or close to being isomorphic. Working with hypersurfaces and emphasizing  examples makes it possible to discuss these questions without any  previous knowledge of algebraic geometry.  At the end we formulate the main recent results and state the most important open questions.

\end{abstract}

  \maketitle

\tableofcontents

Algebraic geometry started as the study of plane curves $C\subset \r^2$ defined  by
a polynomial equation and later extended to surfaces and higher dimensional sets defined by systems of polynomial equations. Besides using $\r^n$, 
it is frequently more advantageous to work with $\c^n$ or with
the corresponding projective spaces $\r\p^n$ and $\c\p^n$.

Later it was realized that the theory also  works if we replace $\r$ or $\c$ by other fields, for example the field of rational numbers $\q$ or even finite fields $\f_q$. When we try to emphasize that the choice of the field is pretty arbitrary, we use  $\a^n$ to denote {\it affine} $n$-space and
$\p^n$ to denote {\it projective} $n$-space. 

Conceptually the simplest algebraic sets are {\it hypersurfaces;} these are defined by 1 equation. That is, 
an {\it affine algebraic hypersurface} of dimension $n$ is the zero set of a 
polynomial $h(x_1,\dots, x_{n+1})$:
$$
X=X(h):=\bigl(h(x_1,\dots, x_{n+1})=0\bigr)\subset \a^{n+1}.
$$
We say that $X$ is {\it irreducible} if $h$ is  irreducible and call $\deg h$ the {\it degree} of $X$.

Working over $\c$,  one can (almost) harmlessly identify the algebraic hypersurface $X$ with the corresponding subset of $\c^{n+1}$,
but already over $\r$ we have to be more careful. For example, 
 $x_1^2+x_2^2+1$  has no real zeros. Thus when we say
``algebraic hypersurface over $\r$'', we  think of not just the real zero set, but also the  complex zero set.  (One could also think about zero sets in even larger fields, for example $\c(t_1,\dots, t_m)$, but it turns out that knowing the zero set over $\c$, or over any algebraically closed field, determines the zero sets over any larger field.) 

Thus in practice we usually think of an affine  hypersurface  $X=X(h)$ as a subset of $\c^{n+1}$ and we keep in mind that it may have been defined by an equation whose coefficients are real or rational. We use $X(\c)$ to denote the set of complex solutions of the equation $h({\mathbf x})=0$.
If $h$ has real or rational coefficients, then it also makes sense to
ask about the real solutions $X(\r)$ or rational solutions $X(\q)$.

We say that two  hypersurfaces are {\it linearly isomorphic} if they differ only  by a   linear change of coordinates. 
In principle one can always decide whether two  hypersurfaces are linearly isomorphic or not, though in practice the computations may be unfeasably long. 

The main  question we   address in these notes is whether there are other, non-linear maps and isomorphisms  between algebraic hypersurfaces. 
The problem naturally divides into two topics. 

\begin{ques} \label{ques.1}
When are two algebraic hypersurfaces $X_1$ and $X_2$ isomorphic?
\end{ques}

\begin{ques} \label{ques.2}
When are two algebraic hypersurfaces $X_1$ and $X_2$ ``almost'' isomorphic?
\end{ques}

Note that we have not yet defined what  an ``isomorphism'' is and 
gave yet no hint what an ``almost'' isomorphism should be. 
Before we give these in Definition~\ref{bir.iso.defn}, we start with an example in Section \ref{sec.1}.
\medskip

{\it Comment on varieties.} Algebraic geometers almost always work with systems of equations, leading to the notions of algebraic sets and varieties. 
For most of the results that we discuss, these form the  natural setting; see \cite{shaf, MR0453732, MR982494, clo-book, dolg-cl, eh-3264} for introductions. 

However, many definitions and results are much easier to state for hypersurfaces, and, once  the foundational material is well  understood,
the basic difficulties are usually very similar.

\begin{ack} These notes are based on lectures given at Northwestern University in April and May, 2018. The hospitality of Northwestern University, and especially of M.~Popa,  gave an ideal time to write them up in expanded form. I received many helpful comments and references from 
 V.~Cheltsov, A.~Corti, T.~de~Fernex, L.~Ein  and Y.~Liu.
 
Partial  financial support    was provided  by  the NSF under grant number
 DMS-1362960 and by the Nemmers Prize of Northwestern University.
\end{ack}

\section{Stereographic projection}\label{sec.1}

The oldest example of an ``almost'' isomorphism between hypersurfaces is the stereographic projection.
It may well have been invented---originally under the name
{\it planisphaerium---}to create starcharts that represent the celestial sphere in a plane.   It has been used in mapmaking  since the XVIth century.

A theorem that  Ptolemy attributes to Hipparchus  ($\sim$ 190-120 BC)
says that stereographic projection maps circles to circles.  
Halley (best known for his comet) proved in  \cite{halley-1695}  that stereographic projection also preserves angles. 
Now we refer to  these properties  by saying that  stereographic projection is conformal.

Instead of these metric properties, we are more interested in the actual formulas that define it.

\begin{say}[Stereographic projection]\label{ster.proj.say}
The formulas for stereographic projection  are nicest if we project
the unit sphere
$$
\s^n:=\bigl(x_1^2+\cdots+x_{n+1}^2=1\bigr)\subset \r^{n+1}
\eqno{(\ref{ster.proj.say}.1)}
$$
from the south pole  ${\mathbf p}_0:=(0,\dots, 0,-1)$ to the plane $x_{n+1}=0$ where we use
$y_1,\dots, y_n$ as coordinates instead of $x_1,\dots,x_n$, 

Geometrically,   pick any point ${\mathbf p}\in \s^n$ (other then the south pole) and
let $\pi({\mathbf p})$ denote the intersection point of the line
$\langle {\mathbf p}_0, {\mathbf p}\rangle$ with the hyperplane $(x_{n+1}=0)$.

Algebraically, we get that 
$$
\pi(x_1,\dots,x_{n+1})=\left(
\frac{x_1}{1+x_{n+1}},\dots, \frac{x_n}{1+x_{n+1}}\right),
\eqno{(\ref{ster.proj.say}.2)}
$$
with inverse
$$
\pi^{-1}(y_1,\dots,y_n)=
\left(
\frac{2y_1}{1+\Sigma},\dots, \frac{2y_n}{1+\Sigma},\frac{1-\Sigma}{1+\Sigma} \right),
\eqno{(\ref{ster.proj.say}.3)}
$$
where $\Sigma=y_1^2+\cdots+y_n^2$. 

Here $\pi$ is not defined at the points where $x_{n+1}=-1$ and
$\pi^{-1}$ is not defined at the points where $y_1^2+\cdots+y_n^2=-1$.
In the real case we get a one-to-one map 
$\s^{n}\setminus(\mbox{south pole})\cong \r^n$  and both
$\pi$  and $\pi^{-1}$  are given by rational functions. 

However, we know that we should also look at the complex case. Note that 
$$
\s^n(\c):=\bigl(x_1^2+\cdots+x_{n+1}^2=1\bigr)\subset \c^{n+1}
$$
is not a sphere,  it is not even compact.
Both indeterminacy sets 
$$
\s^n(\c)\cap( x_{n+1}=-1)\qtq{and} (y_1^2+\cdots+y_n^2=-1)
$$
are $(n-1)$-dimensional hypersurfaces. We get a one-to-one map
between the sets
$$
\s^n(\c)\setminus( x_{n+1}=-1)\qtq{and} \c^n\setminus(y_1^2+\cdots+y_n^2=-1),
$$
which is given by rational functions in both direction.

This will be our definition of a {\it birational map,} which is the
usual name for the  ``almost'' isomorphism of Question \ref{ques.2}.   
\end{say}

A number theoretic variant of stereographic projection is the following.

\begin{say}[A diophantine equation]\label{ster.proj.arith.say}
Consider the diophantine equation
$$
a_1x_1^2+\cdots+a_nx_n^2+x_{n+1}^2=1\qtq{where} a_i\in \q.
\eqno{(\ref{ster.proj.arith.say}.1)}
$$
It determines the quadric hypersurface
$$
Q^n:=\bigl(a_1x_1^2+\cdots+a_nx_n^2+x_{n+1}^2=1\bigr).
\eqno{(\ref{ster.proj.arith.say}.2)}
$$
Again projecting from the south pole  ${\mathbf p}_0:=(0,\dots, 0,-1)$ we get the formulas 
$$
\pi(x_1,\dots,x_{n+1})=\left(
\frac{x_1}{1+x_{n+1}},\dots, \frac{x_n}{1+x_{n+1}}\right)
\eqno{(\ref{ster.proj.arith.say}.3)}
$$
with inverse
$$
\pi^{-1}(y_1,\dots,y_n)=
\left(
\frac{2y_1}{1+\Sigma},\dots, \frac{2y_n}{1+\Sigma},\frac{1-\Sigma}{1+\Sigma} \right),
\eqno{(\ref{ster.proj.arith.say}.4)}
$$
where now $\Sigma=a_1y_1^2+\cdots+a_ny_n^2$. 

Thus we conclude that all rational solutions 
 of the equation
(\ref{ster.proj.arith.say}.1) satisfying  $x_{n+1}\neq -1$ can be written as in (\ref{ster.proj.arith.say}.4) for rational numbers $y_1,\dots, y_n$.

Note that in (\ref{ster.proj.arith.say}.1) we are working with the special case when $a_{n+1}=1$. We could allow instead $a_{n+1}$ to be arbitrary, but then we would need to project from the point $\bigl(0,\dots, 0,-\sqrt{a_{n+1}}\bigr)$. It is not  possible to keep track of the rational solutions this way. 
\end{say}

\section{Projective hypersurfaces}\label{sec.2}

An affine hypersurface of dimension $\geq 1$ is never compact over $\c$, thus we frequently work with its closure in $\c\p^{n+1}$.

\begin{defn}[Projective space]\label{proj.sp.defn}
 We will think of the points of 
$\c\p^{n+1}$ as the lines through the origin in $\c^{n+2}$. 
Thus we use {\it homogeneous coordinates}
$(X_0{:}\cdots{:}X_{n+1})$, where at least one of the $X_i$ most be nonzero. Furthermore, 
$$
(\lambda X_0{:}\cdots{:} \lambda X_{n+1})=(X_0{:}\cdots{:} X_{n+1})\qtq{for all}\lambda \in \c^*.
\eqno{(\ref{proj.sp.defn}.1)}
$$
A polynomial $h(X_0,\dots, X_{n+1})$ can not be evaluated at a point of $\c\p^{n+1}$ since  usually $ h(\lambda X_0{:}\cdots{:} \lambda X_{n+1})\neq h(X_0{:}\cdots{:} X_{n+1})$. However, if  $H$ is homogeneous of degree $d$ then
$$
H(\lambda X_0{:}\cdots{:} \lambda X_{n+1})=\lambda^d H(X_0{:}\cdots{:} X_{n+1}),
\eqno{(\ref{proj.sp.defn}.2)}
$$
thus the zero set of a   homogeneous polynomial is well defined. 
Thus we get {\it projective hypersurfaces}  
$$
X=X(H):=\bigl(H(X_0{:}\cdots{:} X_{n+1})=0\bigr)\subset \p^{n+1}.
\eqno{(\ref{proj.sp.defn}.3)}
$$
One can go between the affine and the projective versions  by the formulas
$$
\begin{array}{rcl}
h(x_1,\dots, x_{n+1})&=&H(1,x_1, \dots, x_{n+1})\qtq{and}\\
H(X_0,\dots, X_{n+1})&=&X_0^{\deg h}\cdot h\bigl(\tfrac{X_1}{X_0},\dots, \tfrac{X_{n+1}}{X_0}\bigr).
\end{array}
\eqno{(\ref{proj.sp.defn}.4)}
$$
Of course any of the $X_i$ could play the special role of $X_0$ above, so
we usually think of $\c\p^{n+1}$ as being covered by  $n+2$ charts, each isomorphic to $\c^{n+1}$.

The {\it cone} over $X$ is the affine hypersurface
$$
C_X=C_{X(H)}:=\bigl(H(X_0,\dots, X_{n+1})=0\bigr)\subset \a^{n+2}.
\eqno{(\ref{proj.sp.defn}.5)}
$$
This gives another way to go between  projective and affine questions. 

({\it Note.} It would be quite convenient to keep the notational distinction between affine coordinates $x_i$ and projective coordinates $X_i$, but people usually use lower case $x_i$ to denote both affine and projective coordinates.) 
\end{defn}

Modern algebraic geometry usually considers the projective variant the basic object and the affine versions as the  local charts on the projective hypersurface.
Note, however, that in algebraic geometry the local charts are very big, they are always dense in the corresponding projective hypersurface.

 \begin{defn} \label{smooth.defn}
A point ${\mathbf p}\in X(G)$ on a projective  hypersurface 
is called {\it smooth} if $\tfrac{\partial G}{\partial x_i}({\mathbf p})\neq 0$ for some $i$. If we are over $\c$ or $\r$, the implicit function therem tells us that 
$X$ is an $n$-dimensional submanifold of  $\c\p^{n+1}$ or of $\r\p^{n+1}$ at its smooth points. Here of course
over $\c$ we count complex dimension, which is twice the real dimension.

The tangent plane  $T_{\mathbf p}X$ at a smooth point
${\mathbf p}:=(p_0{:}\cdots{:}p_{n+1})$ is given by the equation
$$
\tsum_i \tfrac{\partial G}{\partial x_i}({\mathbf p})\cdot x_i=0.
\eqno{(\ref{smooth.defn}.1)}
$$
(Our first inclination would be write $\tsum_i \tfrac{\partial G}{\partial x_i}({\mathbf p})\cdot (x_i-p_i)=0$ instead. This does not seem to make sense as an equation since it is not homogeneous. Luckily, 
$$
\tsum_i \tfrac{\partial G}{\partial x_i}({\mathbf p})\cdot p_i=\deg G \cdot G({\mathbf p})=0
\eqno{(\ref{smooth.defn}.2)}
$$
 since $G$ is homogeneous, thus, after all, we do get (\ref{smooth.defn}.1) this way.)

The other points are  {\it singular.}  Thus the set of all singular points 
is defined by the equations
$$
\tfrac{\partial G}{\partial x_0}=\cdots =\tfrac{\partial G}{\partial x_{n+1}}=0.
\eqno{(\ref{smooth.defn}.3)}
$$ 
It is thus a closed set which is easily seen to be nowehere dense if $G$ is irreducible.
\end{defn}

Over an arbitrary field, we have the following
rather straightforward generalization of the  formulas (\ref{ster.proj.say}.2--3)  for the stereographic projection.

\begin{thm} \label{stere.gen.quad.thm}
Consider the irreducible quadric hypersurface 
$$
Q^n:=\bigl(G(x_0,\dots, x_{n+1})=0\bigr)\subset  \p^{n+1}.
$$
Let ${\mathbf p}\in Q^n$ be a smooth  point with tangent plane $T_{\mathbf p}Q$. Pick any hyperlane
$H\cong \p^n$ that does not contain ${\mathbf p}$. Then the following hold.
\begin{enumerate}
\item  Projection of $\p^{n+1}$ from ${\mathbf p}$ to $H$ gives a one-to-one map
$$
\pi:  Q^n\setminus T_{\mathbf p}Q  \longrightarrow  H\setminus T_{\mathbf p}Q \cong \a^n.
$$
\item  The coordinate functions of $\pi$  are quotients of linear polynomials.
\item  The coordinate functions of $\pi^{-1}$  are quotients of quadratic polynomials.
\item If the coefficients of $G$ are in a field $k$ then so are the
coefficients of $\pi$ and $\pi^{-1}$. \qed
\end{enumerate}
\end{thm}

\section{Rational and birational maps}\label{sec.3}

\begin{defn}[Rational maps, affine case]\label{rat.map.aff.defn}
The basic functions in algebraic geometry are polynomials
$p(x_1,\dots, x_n)$ and their quotients
$$
\phi(x_1,\dots, x_n)=\tfrac{p(x_1,\dots, x_n)}{q(x_1,\dots, x_n)},
$$
called {\it rational functions.} Since a polynomial ring has unique fatorization, we may assume that $p, q$ are relatively prime. In this case
  $\phi$ is  not  defined  along the hypersurface $(q=0)$.  

A   {\it rational map} from  $\a^{n}$ to
$\a^{m}$ is a map given by  rational coordinate functions
$$
\Phi(x_1,\dots, x_{n})=\bigl(\phi_1(x_1,\dots, x_{n}), \dots , \phi_{m}(x_1,\dots, x_{n})\bigr),
$$
where $\phi_i=p_i/q_i$. As before, $\Phi$ need not be  everywhere defined.
\end{defn}

\begin{defn}[Rational maps, projective case]\label{rat.map.pr.defn}
As we already noted in (\ref{proj.sp.defn}.2), a  rational function  $p(x_0,\dots, x_n)$ can not be evaluated at a point of $\p^n$, but if it
is a quotient of 2 homogeneous polynomials of the same degree then
$$
\phi(\lambda x_0{:}\dots{:} \lambda x_n)=\phi(x_0{:}\cdots{:} x_n),
\eqno{(\ref{rat.map.pr.defn}.1)}
$$
and $\phi$ defines a  {\it rational function} on $\p^n$.

 A {\it rational} map from  $\p^{n}$ to
$\p^{m}$ is a map given by  rational functions
$$
\Phi(x_0{:}\cdots{:} x_{n})=\bigl(\phi_0(x_0{:}\cdots{:} x_{n}):\cdots : \phi_{m}(x_0{:}\cdots{:} x_{n})\bigr).
\eqno{(\ref{rat.map.pr.defn}.2)}
$$
Where is $\Phi$ defined?
Note first that  each $\phi_i$ is a quotient of 2 polynomials  $\phi_i=p_i/q_i$
thus  $\phi_i$ might  not  be defined along the hypersurfaces $(q_i=0)$. 
Furthermore,  every point in $\p^{m}$ has at least 1 non-zero coordinate, so  
$\Phi$ is  also not defined  along the intersection  of the hypersurfaces $(p_0=0)\cap\cdots\cap (p_{n}=0)$. 

However,  the coordinate functions of a rational map are not unique
since
$$
\bigl(\phi_0({\mathbf x}):\cdots : \phi_{m}({\mathbf x})\bigr)=\bigl(\psi({\mathbf x})\phi_0({\mathbf x}):\cdots : \psi({\mathbf x})\phi_{m}({\mathbf x})\bigr)
\eqno{(\ref{rat.map.pr.defn}.3)}
$$
for any rational function $\psi({\mathbf x}) $. This is confusing at the beginning, but
we can also turn the non-uniqueness of projective coordinates to our advantage.
Multiplying through with the least common denominator of the $\phi_i$ in (\ref{rat.map.pr.defn}.2) we represent $\Phi$ as
$$
\Phi({\mathbf x})=\bigl(p_0({\mathbf x}):\cdots : p_{m}({\mathbf x})\bigr),
\eqno{(\ref{rat.map.pr.defn}.4)}
$$
where the $p_i$ are polynomials. 
Using that  polynomial rings have unique factorization, we can cancel out common factors and write
$$
p_i({\mathbf x})=c_i\tprod_j r_j({\mathbf x})^{m(i,j)},
\eqno{(\ref{rat.map.pr.defn}.5)}
$$
where the $c_i$ are constants, the $r_j$ are irreducible polynomials and $m(i,j)\geq 0$. Furthermore, at least one $m(i,j)$ is zero for every $j$.
Then $\Phi$ is defined  outside the common zero set of the $p_i$, which is contained in 
$$
Z:=\cup_{j_1\neq  j_2} \bigl(r_{j_1}({\mathbf x})=r_{j_2}({\mathbf x})=0\bigr).
\eqno{(\ref{rat.map.pr.defn}.6)}
$$
This is made up of many pieces, but each of them is defined by the vanishing of 2 relatively prime polynomials. Thus we expect that $Z$ has codimension $\geq 2$ in $\p^n$. This is not hard to see; we state a
  general version as Theorem~\ref{rtl.cod.2.thm}.
\end{defn}

\begin{defn}[Rational maps of hypersurfaces]\label{rat.map.hyp.defn}
A {\it rational function} on a hypersurface $X\subset \p^{n}$ is the restriction of a rational function $\phi$ on $\p^n$  to $X$, provided the restriction makes sense. That is, when $\phi$ is defined
on a dense open subset of $X$.

As before,  a {\it rational map} $\Phi$ from a hypersurface $X\subset \p^{n}$ to
$\p^{m}$ is given by rational functions  $(\phi_0{:}\cdots{:}\phi_m)$. 
We denote rational maps by a dashed arrow $\map$. 
If $\Phi(X)$ is contained in a hypersurface $Y\subset \p^m$, then
$\Phi$ defines a rational map  $\Phi:X\map Y$. 
(If we work over $\c$, then by  $\Phi(X)\subset Y$ we mean that
$\Phi({\mathbf p})\in Y$ whenever $\Phi$ is defined at ${\mathbf p}\in X$. The real version of this  is not   the correct definition  since 
$X(\r)$ may be empty. Thus we need to work with complex points or use the algebraic version given in Remark~\ref{aff.bir.rem}.)

We use  3 important properties of rational functions.
\medskip

{\it \ref{rat.map.hyp.defn}.1.}
 If $X$ is irreducible and $\phi$ is defined at a single point of $X$ then it is defined
on a dense open subset of $X$, cf.\ \cite[Sec.I.3.2]{shaf}. 

\medskip

{\it \ref{rat.map.hyp.defn}.2.} If $X\subset \a^{n+1}$ is affine and $\phi|_X$ is everywhere defined then there is a polynomial $p$ on $\a^{n+1}$ such that
$p|_X=\phi|_X$.   For example,  $\frac{x^3-y}{x+y}$ has poles but its restriction to the $x$-axis is everywhere defined. We can take $p:=x^2$.

This is usually proved as a consequence of Hilbert's Nullstellensatz, see \cite[Sec.I.3.2]{shaf}. 

\medskip

{\it Purity Principle \ref{rat.map.hyp.defn}.3.}  Let $\phi$ be a rational function on a hypersurface $X$. Then the set of points where $\phi$ is not defined has codimension 1 everywhere.

This is more subtle, in books it  is usually treated as a combination of
Krull's principal ideal theorem and Serre's $S_2$ property for hypersurfaces.
However this form goes back to Macaulay \cite{macaulay-1916}. 
 It is an algebraic counterpart of Hartogs's extension theorem  in complex analysis.
\medskip

\end{defn}

\begin{defn}[Morphisms]  \label{rat.map.def.defn}
$\Phi$ is  {\it defined} at a point ${\mathbf p}\in X$ if it has some representation  
$$
\Phi({\mathbf x})=\bigl(\psi({\mathbf x})\phi_0({\mathbf x}):\cdots : \psi({\mathbf x})\phi_{m}({\mathbf x})\bigr)
\eqno{(\ref{rat.map.def.defn}.1)}
$$
as in (\ref{rat.map.pr.defn}.3), where every coordinate function is defined at ${\mathbf p}$ and at least one of them is nonzero.
(The traditional terminology is  ``$\Phi$ is  regular'' at ${\mathbf p}$,  but be warned that
``regular''  has other, conflicting uses, even in algebraic geometry.)

A {\it morphism} is a rational map that is everywhere defined.
We denote morphisms by a solid arrow $\to$. 

\medskip

{\it Warning.} Unlike for $\p^n$, where we could write down the optimal
representation in (\ref{rat.map.pr.defn}.5), there need not be a single
form (\ref{rat.map.def.defn}.1) that shows regularity at every point.

For example, consider the affine hypersurface $X=(x_1x_2=x_3x_4)\subset \a^4$
and the rational function  $\phi=x_1/x_3$. It can have poles only along
$x_3=0$. Next note that $\phi=x_4/x_2$ and the latter form can have poles only along
$x_2=0$. So in fact we have poles only along the 2-plane  $(x_2=x_3=0)$.

\end{defn}

Nonetheless, the following generalization of (\ref{rat.map.pr.defn}.6) still holds,
see \cite[Sec.II.3.1]{shaf}.

\medskip

\begin{thm}\label{rtl.cod.2.thm}
 Let $X$ be a smooth hypersurface  and
$\Phi:X\map \p^m$ a rational map. Then there is a closed subset $Z\subset X$ of codimension $\geq 2$  such that $\Phi$ is defined on $X\setminus Z$. \qed
\end{thm}

Now we come to the definition of ``isomorphism'' and ``almost isomorphism.''

\begin{defn}\label{bir.iso.defn} Let $\Phi:X\map Y$ be  a rational  map between hypersurfaces.
 \begin{enumerate}
\item $\Phi$  is   {\it birational} if there is a rational map $\Phi^{-1}:Y\map X$ that is the inverse of $\Phi$.  That is, $\Phi^{-1}\circ \Phi$ and $\Phi\circ \Phi^{-1}$ are both the identity, wherever they are defined.
If this holds then we say that $X$ and $Y$ are  {\it birationally equivalent} or  {\it birational.}
\item An $n$-dimensional  hypersurface $X$ is called {\it rational} 
if it is birational to $\p^n$. (In the  old  literature this is frequently called birational.)
\item $\Phi$  is  an 
{\it isomorphism} if both $\Phi$ and $\Phi^{-1}$ are morphisms, that is, everywhere defined. 
\end{enumerate}
\end{defn}

\begin{rem}\label{aff.bir.rem}
Note that an affine hypersurface is birational to its projective closure, thus birationality can be checked using affine equations. Two  hypersurfaces 
$$
\begin{array}{rcl}
X&:=&\bigl(h(x_1,\dots, x_{n+1})=0\bigr)\subset \a_{\mathbf x}^{n+1}\qtq{and}\\
Y&:=&\bigl(g(y_1,\dots, y_{n+1})=0\bigr)\subset \a_{\mathbf y}^{n+1}
\end{array}
$$
are  birational  if there are rational maps
$$
\begin{array}{l}
\Phi=\bigl(\phi_1(x_1,\dots, x_{n+1}),\dots, \phi_{n+1}(x_1,\dots, x_{n+1})\bigr)\qtq{and}\\
\Psi=\bigl(\psi_1(y_1,\dots, y_{n+1}),\dots, \psi_{n+1}(y_1,\dots, y_{n+1})\bigr)
\end{array}
$$
with the following properties.

(\ref{aff.bir.rem}.1) $\Phi$ maps $X$ to $Y$. In terms of equations this means that
$g\bigl(\phi_1,\dots, \phi_{n+1}\bigr)$
vanishes on $X$. That is, if we write 
$g\bigl(\phi_1,\dots, \phi_{n+1}\bigr)$ as the quotient of 2 relatively prime polynomials then the 
 numerator is divisible by $h$.

(\ref{aff.bir.rem}.2) $\Psi$ maps $Y$ to $X$.

(\ref{aff.bir.rem}.3) 
$\Phi: X\map  Y$ and  $\Psi: Y\map X$
are inverses if each other.  Note that 
$$
\Psi\circ \Phi=
\bigl(\psi_1(\phi_1,\dots, \phi_{n+1}),\dots, \psi_{n+1}(\phi_1,\dots, \phi_{n+1})\bigr),
$$
and this is the identity on $X$ iff $1-\psi_i(\phi_1,\dots, \phi_{n+1})$ vanishes on $X$ for every $i$. Similarly,  $1-\phi_j(\psi_1,\dots, \psi_{n+1})$ vanishes on $Y$ for every $j$. 
\end{rem}

\section{The main questions}\label{sec.4}

Now we are ready to formulate the precise versions of
Questions~\ref{ques.1}--\ref{ques.2}.
In both cases the ideal complete answer would consists of 2 steps.
\medskip

$\bullet$ Describe a set of ``elementary'' isomorphisms/birational maps between hypersurfaces.

$\bullet$ Prove that every isomorphism/birational map between hypersurfaces is ``elementary,'' or at least  a composite of  ``elementary'' maps.
\medskip

So what are these  ``elementary'' isomorphisms/birational maps?

\medskip
{\bf Isomorphisms of projective hypersurfaces}
\medskip

For projective hypersurfaces $X\subset \p^{n+1}$ the ``elementary'' isomorphisms are those that  are induced by an automorphism of $ \p^{n+1}$. It is not hard to see that  $ \aut(\p^{n+1})\cong \PGL_{n+2}$ (cf.\ \cite[III.1.Exrc.17]{shaf}), so these are exactly the 
linear isomorphisms. 
This leads to  the projective version of Question~\ref{ques.1}

\begin{ques} \label{ques.1.1} Let $X_1, X_2$ be projective algebraic hypersurfaces in $\p^{n+1}$ and $\Phi:X_1\to X_2$ an isomorphism. Is $\Phi$ linear?
That is, is there an automorphism $\Psi\in \aut(\p^{n+1})\cong \PGL_{n+2}$
such that $\Phi=\Psi|_{X_1}$?
\end{ques}

We discuss almost complete answers to this is Section \ref{sec.6}. 

\medskip
{\bf Isomorphisms of affine hypersurfaces}
\medskip

The projective case  suggests that for affine  hypersurfaces the ``elementary'' isomorphisms  should be  those that are induced by an automorphism of $\a^{n+1}$.
A major difficulty  is that $\aut(\a^{n+1})$ is  infinite dimensional. For example, given any polynomials
$g_i(x_{i+1},\dots, x_{n+1})$, the map
$$
\Phi:(x_1,\dots, x_{n+1})\mapsto
\bigl(x_1+g_1,x_2+g_2, \dots, x_{n+1}+g_{n+1})
$$
is an  automorphism of $\a^{n+1}$. These automorhisms and $\GL_n$ together generate the {\it tame subgroup} of $\aut(\a^{n+1})$. If $n=1$ then we get the whole $\aut(\a^2)$ by \cite{MR542446} but  the analogous result does not hold in higher dimensions. A counter example was proposed by Nagata  \cite{MR0337962} and proved in  \cite{MR2015334}; see \cite{MR3561400} for a survey and  \cite{lempert} for a complex analytic version.

Already the simplest case of the affine isomorphism problem is a quite nontrivial result of \cite{MR0379502} and \cite{MR0338423}.

\begin{thm}[Abhyankar-Moh-Suzuki] Let $C\subset \a^2$ be a curve
and $\Phi:C\to \a^1$ an isomorphism. Then $\Phi$ extends to an isomorphism
$$
\Psi: \bigl(C\subset \a^2\bigr)\cong \bigl((\mbox{\rm coordinate axis})\subset \a^2\bigr). \qed
$$
\end{thm}

The higher dimensional generalization of this is the Abhyankar-Sathaye Conjecture, about which very little is known and  a series of counter examples is proposed in \cite{MR1282208}; see  also \cite[Chap.5]{MR1790619}.

\begin{conj} Let $X\subset \a^{n+1}$ be a hypersurface
and $\Phi:X\to \a^n$ an isomorphism. Then $\Phi$ extends to an isomorphism
$$
\Psi: \bigl(X\subset \a^{n+1}\bigr)\cong \bigl((\mbox{\rm coordinate hyperplane})\subset \a^{n+1}\bigr). 
$$
\end{conj}

Interestingly, the union of the $n$ coordinate hyperplanes 
$(x_1\cdots x_n=0)$ 
has a unique embedding into  $\c^n$ by
\cite{MR1446200}. 
For some recent results connecting the  affine isomorphism problem
with the  methods  of Section~\ref{sec.10},  see \cite{cdp2017}. 

\medskip
{\bf Birational equivalence of hypersurfaces}
\medskip

It seems natural to formulate Question~\ref{ques.2} in two variants.
Note that an affine hypersurface is birational to its projective closure, so
for the birationality questions there is  no need to distinguish the affine and the projective cases.

\begin{ques} \label{ques.2.1}
Which algebraic hypersurfaces are rational?
\end{ques}

\begin{ques} \label{ques.2.2}
When are two algebraic hypersurfaces birational to each other?
\end{ques}

\begin{exmp}[Rational hypersurfaces] \label{rtl.sing.exmps}
It has been long understood that even high degree hypersurfaces can be rational if they are very singular. Here are some examples of this.
\smallskip

(\ref{rtl.sing.exmps}.1) Let $X\subset \p^{n+1}$ be given by an equation
$$
H_{d-1}(x_0,\dots, x_n)x_{n+1}+ H_{d}(x_0,\dots, x_n)=0
$$
where $H_j$ is homogeneous of degree $j$. Projection from the point
$(0{:}\cdots{:}0{:}1)$ to the hyperplane $(x_{n+1}=0)$ is birational, the inverse is given by
$$
(x_0{:}\cdots{:} x_n)\mapsto \Bigl(x_0{:}\cdots{:} x_n{:}\tfrac{H_{d}(x_0,\dots, x_n)}{H_{d-1}(x_0,\dots, x_n)}\Bigr).
$$
If $d=2$ then this is the stereographic projection. If $d\geq 3$ then
$(0{:}\cdots{:}0{:}1)$ is a singular point of $X$. 
\smallskip

(\ref{rtl.sing.exmps}.2) 
 Let $X\subset \p^{2n+1}$ be a  hypersurface of degree $2d+1$ given by an equation of the form
$$
\tsum_{I,J,k}\ a_{I,J,k}\cdot M_I\cdot N_J\cdot x_k,
$$  where $M_I$ is a degree $d$ monomial in the variables
$x_0,\dots, x_{n}$ and $N_J$ is a degree $d$ monomial in the variables
$x_{n+1},\dots, x_{2n+1}$ and $x_k$ is arbitrary. 
(In particular,  $X$  contains the  linear subspaces 
$L_1:=(x_0=\cdots= x_{n}=0)$ and $L_2:=(x_{n+1}=\cdots= x_{2n+1}=0)$.) Then $X$ is rational. We will work out the $2d+1=3$ case in detail in Proposition~\ref{cubic.2L.rtl.prop}.

If $2d+1=3$ then a general such $X$ is smooth, but for $2d+1\geq 5$ it is always singular along the $L_i$.

\smallskip

(\ref{rtl.sing.exmps}.3) Take  degree $d$ homogeneous polynomials
$p_0({\mathbf x}), \dots,  p_{n+1}({\mathbf x})$. 
They define a rational map $\Phi: \p^n\map \p^{n+1}$. 
It is not hard to show that if the $p_i$ are general then
$\Phi$ is a morphism whose image $Y=Y(p_0,\dots, p_{n+1})$ is a rational hypersurface of degree
$d^n$. However, $Y$ has very complicated self-intersections.
\end{exmp}

These and many other examples suggest that
Questions~\ref{ques.2.1}--\ref{ques.2.2} are reasonable only if
the hypersurfaces are smooth or mildly singular. We will focus on the smooth cases.
 
\begin{ques} \label{ques.2.1.s}
Which smooth algebraic hypersurfaces are rational?
\end{ques}

\begin{ques} \label{ques.2.2.s}
When are two smooth algebraic hypersurfaces birational to each other?
\end{ques}

So far is not even clear that not all hypersurfaces are rational.
We see later that this is not at all the case, but first  we give  some rational examples.

\section{Rationality of cubic hypersurfaces}\label{sec.5}

We proved in Theorem~\ref{stere.gen.quad.thm} that a smooth quadric hypersurface is rational over a  field $k$ 
if it has at least one  $k$-point that we can project it from.  
In this section we consider the rationality of  cubic hypersurfaces.
This is a topic with many interesting results and still unsolved questions.
The simplest rationality construction  is the following special case of (\ref{rtl.sing.exmps}.2). 

\begin{prop}\label{cubic.2L.rtl.prop}
 Let $S\subset \p^3$ be a cubic surface that contains 2 disjoint lines. Then $S$ is rational.
\end{prop}

Geometric proof.  Let the lines be $L_1, L_2$. Pick points ${\mathbf p}_i\in L_i$.
The line connecting them  $\langle {\mathbf p}_1, {\mathbf p}_2\rangle$ meets $S$ in 3 points.
We already know 2 of them, namely ${\mathbf p}_1, {\mathbf p}_2$.  Let $\phi({\mathbf p}_1, {\mathbf p}_2)$ be the third intersection point. This gives a rational map
$$
\phi: L_1\times L_2\map S.
\eqno{(\ref{cubic.2L.rtl.prop}.1)}
$$
To get its inverse, pick ${\mathbf q}\in \p^3$ and let $\pi:\p^3\map \p^2$  denote the projection from ${\mathbf q}$. Then $\pi(L_1), \pi(L_2)$ are 2 lines in $\p^2$, hence they meet at a unique point ${\mathbf q}'$. Thus we get ${\mathbf p}_i\in L_i$ such that
$\pi({\mathbf p}_i)={\mathbf q}'$. This gives
 $$
\psi: \p^3 \map L_1\times L_2 \qtq{such that} \psi|_S=\phi^{-1}. \qed
\eqno{(\ref{cubic.2L.rtl.prop}.2)}
$$

Algebraic proof. We can choose coordinates such that
$L_1=(x_0=x_1=0)$ and $L_2=(x_2=x_3=0)$. Thus the equation of $S$ can be  written (non-uniquely) as
$$
\tsum_{i=0,1}\tsum_{j=2,3} \ell_{ij}({\mathbf x})x_ix_j,
\eqno{(\ref{cubic.2L.rtl.prop}.3)}
$$
where the $\ell_{ij}({\mathbf x})$ are linear in the variables $x_0,\dots, x_3$.

If $p_1=(a_0{:}a_1{:}0{:}0)$ and $p_2=(0{:}0{:}a_2{:}a_3)$
then the line connecting them is 
$$
\p^1\ni (s{:}t)\mapsto (sa_0{:}sa_1{:}ta_2{:}ta_3) \in \p^3.
\eqno{(\ref{cubic.2L.rtl.prop}.4)}
$$
Substituting  into (\ref{cubic.2L.rtl.prop}.3) we get an equation
$$ \underset{i=0,1}{\tsum}\underset{j=2,3}{\tsum}\ell_{ij}(sa_0{:}sa_1{:}ta_2{:}ta_3)sa_ita_j=
st  \underset{i=0,1}{\tsum}\underset{j=2,3}{\tsum}\ell_{ij}(sa_0{:}sa_1{:}ta_2{:}ta_3)a_ia_j.
\eqno{(\ref{cubic.2L.rtl.prop}.5)}
$$
After dividing by $st$, the remaining equation is linear in $s,t$, thus it has a unique solution  (up to a multiplicative constant)
$$
(s{:}t)=
\Bigl(
- \underset{i=0,1}{\tsum}\underset{j=2,3}{\tsum}\ell_{ij}(0{:}0{:}a_2{:}a_3)a_ia_j \ :\  \underset{i=0,1}{\tsum}\underset{j=2,3}{\tsum}\ell_{ij}(a_0{:}a_1{:}0{:}0)a_ia_j\Bigr). \qed
\eqno{(\ref{cubic.2L.rtl.prop}.6)}
$$

This can be used to prove the following result of 
\cite{MR1579329}.

\begin{cor}[Clebsch] \label{cubic.C.rtl.prop}
Every smooth cubic surface over $\c$ is rational.
\end{cor}

Proof. By Proposition~\ref{cubic.2L.rtl.prop} it is enough to find 2 disjoint lines. In fact, every
smooth cubic surface over $\c$ contains 27 lines. This was first proved by
Cayley and Salmon in 1849 and published in \cite{salmon-2-old}; see the next example or any of \cite[Chap.7]{MR982494}, \cite[Chap.1]{MR1442522} and
\cite[Sec.IV.2.5]{shaf} for  proofs.\qed

\begin{exmp} \label{lines.on.diag.exmp}
Consider  the degree $d$ hypersurface
$$
X:=\bigl(x_0^d+\cdots + x_{2n+1}^d=0\bigr)\subset \p^{2n+1}.
$$
Divide the indices into $n+1$ disjoint ordered pairs
$(x_{\tau(i)}, x_{\sigma(i)}):i=0,\dots, n$. 
Fix  $d$th roots 
$\epsilon_0,\dots, \epsilon_n$ of $-1$. Then
$$
L=L(\tau, \sigma, \epsilon):=\bigl(x_{\tau(i)}=\epsilon_ix_{\sigma(i)}:i=0,\dots, n\bigr)
$$
is a linear space of dimension $n$ contained in $X$.
Thus $X$ contains at least 
$$
\frac{d^{n+1}(2n+2)!}{2^{n+1}(n+1)!}
$$ linear subspaces of dimension $n$.
Not all but many pairs of these  linear subspaces are disjoint, for example
$L(\tau, \sigma, \epsilon)$ is disjoint from $L(\tau', \sigma', \epsilon')$
if $\tau=\tau', \sigma=\sigma'$ but $\epsilon_i\neq  \epsilon'_i$ 
for every $i$.
For cubic surfaces we get 27 lines. 
\end{exmp}

\medskip

The following result of \cite{MR0009471}  is especially strong over number fields; see \cite[Chap.2]{ksc} for a modern treatment.

\begin{thm}[Segre] The cubic surface
$$
S=S(a_0,a_1,a_2,a_3):=\bigl(a_0x_0^3+a_1x_1^3+a_2x_2^3+a_3x_3^3=0\bigr)\subset \p^3
$$
defined over a field $k$ 
is not rational over $k$ if for every  permutation $\sigma$ of the indices,  the quotient
$(a_{\sigma(0)}a_{\sigma(1)})/(a_{\sigma(2)}a_{\sigma(3)})$ is not a
 cube in $k$. \qed
\end{thm}

\begin{say}[Cubic 3-folds] In  contrast to
Corollary~\ref{cubic.C.rtl.prop}, 
   smooth cubic 3-folds $X_3\subset \p^4$ are all   non-rational \cite{CG72}.
This is proved by analyzing the Hodge structure of the cohomology groups.
This method has been very successful for many other 3-dimensional cases, but extensions to higher dimensions are still lacking.
\end{say}

\begin{say}[Cubic 4-folds] A cubic 4-fold is given by a degree 3 homogeneous polynomial in 6 variables. These form a vector space of dimension   
$\binom{8}{5}=56$. Since we do not care about multiplicative constants,
we can think of the set of all cubic 4-folds as points in a $\p^{55}$. 
The smooth ones correspond to an open subset of it. 
(This is a general feature of algebraic geometry. Isomorphism classes frequently naturally correspond to points of another algebraic variety, called the {\it moduli space.})

 A long standing conjecture that grew out of the works of
\cite{MR0020273, MR0010433}
says  that 
a general smooth cubic 4-fold  is not rational.  \cite{MR1658216} describes  countably many hypersurfaces $H_i\subset \p^{55}$
and  conjectures  that  a cubic 4-fold $X$ is not rational if the corresponding point $[X]\in \p^{55}$ is outside these hypersurfaces. See the collection 
\cite{MR3587799} for recent surveys.

The strongest results, due to Russo and Staglian\'o \cite{rus-sta}
prove rationality for  cubic 4-folds corresponding to 3 of the hypersurfaces $H_i$. 
Just to show that there are some rather subtle rationality constructions, here is one of the  beautiful examples discovered by them.

Pick 10 general points  $P_j\in \p^2$ and 6 general polynomials
$p_i(x,y,z)$ of degree 10 that vanish  with multiplicity 3 at  all the points $P_j$. Let $F$ be the surface obtained as the image of the map
$$
\Psi: \p^2\to \p^5: (x{:}y{:}z)\mapsto
\bigl(p_0(x{:}y{:}z):\cdots:p_5(x{:}y{:}z)\bigr).
$$
\cite{rus-sta} shows that there are 5 linearly independent quintic polynomials 
$q_0,\dots, q_4$ on $\p^5$ that vanish along $F$ with multiplicity 2. 
These give a rational map
$\Phi: \p^5\map \p^4$ whose general fiber is a conic that intersects  $F$ at 5 points. This implies that if $X$ is a cubic hypersurface that contains $F$ then
we get a birational map
$$
\Phi|_X: X\map \p^4.
$$
It turns out that this construction  works for a 54-dimensional subfamily of the   55-dimensional family of all cubics.

The method of the proof of these claims in \cite{rus-sta} shows both the power and the limitations of algebraic geometry. Consider for instance the claim about the 5 linearly independent quintic polynomials. Standard methods imply that there are at least 5 such linearly independent quintics. Then other semicontinuity arguments show that if there are exactly  5 in at least one concrete example,  then there are exactly  5 in ``almost all'' cases. 
Moreover, it is enough to find one such example over a finite field. We can thus ask a computer to try out random cases until it hits one that works.  
Since in algebraic geometry ``almost all'' usually means a dense open subset, there is every reason to believe that this leads to a solution. 

The disadvantage is that while we proved our claims for almost all cases, we do not yet know for which ones. Answering the latter question may need new ideas or techniques.
\end{say}

\section{Isomorphism of  hypersurfaces}\label{sec.6}

 The following theorem answers Question~\ref{ques.1.1}.

\begin{thm}\label{isom.hyp.thm}
 Let $X_1, X_2\subset \p^{n+1}$ be  irreducible hypersurfaces 
and $\Phi:X_1\cong X_2$ an isomorphism. Then $\Phi$ is obtained by a linear change of coordinates in $\p^{n+1}$, except possible in the following 3 cases.
\begin{enumerate}
\item $\dim X_1=\dim X_2=1$ and $\{\deg X_1, \deg X_2\}=\{1,2\}$, 
\item $\dim X_1=\dim X_2=1$ and $\deg X_1=\deg X_2=3$ or
\item $\dim X_1=\dim X_2=2$ and $\deg X_1=\deg X_2=4$.
\end{enumerate} 
\end{thm}

First we  describe in detail the exceptional cases.

The first exceptional case is given by the stereographic projection
of a plane conic to a line. This is the only dimension where the stereographic projection is an isomorphism. 

The second exceptional case is $d=3$ and $n=1$, that is, degree 3 curves in $\p^2$. This is the theory of elliptic curves and elliptic integrals. 
We give a description using the theory of the Weierstrass $\wp$-function;
see \cite[Chap.1]{MR0257326}, \cite[Chap.6]{silverman},    \cite[Chap.9]{stein-complex} or many other books on complex analysis for details.

\begin{exmp}\label{ell.aut.exmp.1} 
Let $\Lambda\subset \c$ be a lattice; after multiplying by a suitable $c\in \c^*$ we can achieve  that its generators are $1$ and $\tau$ where $\operatorname{Im}(\tau)>0$. The Weierstrass  $\wp(z)$ and its derivative $\wp'(z)$ satisfy an equation
$$
\wp'(z)^2=4 \wp^2(x)-g_4\wp(z)-g_6,
$$
where the precise formulas for the $g_i$ are not important for us.
Thus  $\sigma: z\mapsto \bigl(\wp(z){:} \wp'(z){:}1\bigr)$
gives an isomorphism
$$
\sigma:\c/\Lambda\cong C_{\tau}:=\bigl(y^2w=4x^3-g_4xw^2-g_6w^3)\subset \p^2_{xyw}.
$$
The origin is mapped to the point at infinity which is an  inflection point of the curve.

For any  $z_0\in \c$, we can also use  
$\sigma_{z_0}: z\mapsto \bigl(\wp(z-z_0){:} \wp'(z-z_0){:}1\bigr)$
  to get another isomorphism
$$
\sigma_{z_0}:\c/\Lambda\cong C_{\tau}=\bigl(y^2w=4x^3-g_4xw^2-g_6w^3)\subset \p^2_{xyw}.
$$
Here the point $z_0\in \c$  is mapped to the inflection point at infinity.

A linear automorphism of $\p^2_{xyw}$ preserves inflection points, so
the automorphism $\sigma_{z_0}\circ \sigma^{-1}:C_{\tau}\cong C_{\tau}$ can not be linear 
if  $\sigma(z_0)$ is not an  inflection point. 

With a little more work we get that $\sigma_{z_0}\sigma^{-1}:C_{\tau}\cong C_{\tau}$ is linear iff  $z_0\in \tfrac13 \Lambda$. 
\end{exmp}

\begin{exmp}\label{ell.aut.exmp.2}  A more geometric way of obtaining a non-linear automorphism of a smooth plane cubic $C\subset \p^2$ is the following. Pick a point $p_0\in C$ and for $p\in C$ let $\tau_{p_0}(p)$ denote the 3rd intersection point of the line $\langle p, p_0\rangle$ with $C$.
Thus $\tau_{p_0}$ is an involution of $C$ and it is linear iff $p_0$ is an
inflection point of $C$. 

In the description of Example~\ref{ell.aut.exmp.1}, we can get
 $\tau_{p_0}$ as  $z\mapsto 2z_0-z$ where $\sigma(z_0)=p_0$.
\end{exmp}

The third exceptional case is $d=4$ and $n=2$, that is, degree 4 surfaces in $\p^3$. Here there are many known examples of non-linear isomorphisms but these are not easy to find and we do not have a complete description. A rather large set of examples is the following. 

\begin{exmp} \label{det.K3.nonisom.exmp}
Start with $\p^3_{\mathbf x}\times \p^3_{\mathbf y}$ and
4 bilinear hypersurfaces given by equations
$H_k:=(\sum_{ij} a_{ij}^kx_iy_j=0)$.  Set $S:=H_1\cap\cdots\cap H_4$. 
In order to compute its projection to $\p^3_{\mathbf x}$, we view
the  $H_k$  as linear equations in the ${\mathbf y}$-variables. 
The coefficient matrix is
$$
{\mathbf B}=\bigl(b_{kj}\bigr)\qtq{where}  b_{kj}= \tsum_{i} a_{ij}^kx_i.
$$
A given point  ${\mathbf p}\in\p^3_{\mathbf x}$ is in $\pi_{\mathbf x}(S)$ iff
 ${\mathbf B}({\mathbf p}){\mathbf y}^t=0$ has a nonzero solution ${\mathbf y} $ and this in turn holds iff 
$\det {\mathbf B}({\mathbf p})=0$.
Thus
$$
\pi_{\mathbf x}(S)=\bigl(\det {\mathbf B}({\mathbf x})=0\bigr)\subset \p^3_{\mathbf x}.
$$
Similarly, the image of the projection to $\p^3_{\mathbf y}$ is given by
$$
\pi_{\mathbf y}(S)=\bigl(\det {\mathbf C}({\mathbf y})=0\bigr)
\qtq{where} {\mathbf C}=\bigl(c_{ki}\bigr)= \bigl(\tsum_{j} a_{ij}^ky_j\bigr).
$$
It is not hard to check that, for general $a_{ij}^k$, the projections
give isomorphisms
$$
S\cong \bigl(\det {\mathbf B}({\mathbf x})=0\bigr)
\qtq{and} 
S\cong\bigl(\det {\mathbf C}({\mathbf y})=0\bigr).
$$
In particular, the quartic surfaces 
$\bigl(\det {\mathbf B}=0\bigr) $ and $\bigl(\det {\mathbf C}=0\bigr) $ are  isomorphic to each other.
By Cramer's rule  the isomorphism between them can be  given  by the coordinate functions 
$$
\phi_j({\mathbf x})=(-1)^{j-1}\det {\mathbf B}_{4j}({\mathbf x}), 
$$
where ${\mathbf B}_{4j} $ is the submatrix obtained by 
removing the 4th row and $j$th column. So the  $\phi_j({\mathbf x})$ are  cubic polynomials.  

In principle it could  happen that  $\bigl(\det {\mathbf B}({\mathbf x})=0\bigr)$ and $\bigl(\det {\mathbf C}({\mathbf y})=0\bigr)$ are
isomorphic by a linear isomorphism.  This does not happen for general $a_{ij}^k$, though this needs proof.
See \cite{2016arXiv160406265S, 2016arXiv160204588O} for this and other such examples.

One can do the same game with $\p^{n+1}_{\mathbf x}\times \p^{n+1}_{\mathbf y}$ and
$n+2$ bilinear hypersurfaces. It turns out  that for $n\geq 3$ the projections
are birational but they are not isomorphisms and the image hypersurfaces are   always singular. Nonetheless, this construction gives very interesting examples of mildly singular  (in fact with terminal singularities)  hypersurfaces that are birational to each other in an unexpected way.
\end{exmp}

\begin{say}[Comments on the proof of Theorem~\ref{isom.hyp.thm}]\label{isom.hyp.thm.pf.1}
The theorem  is an easy consequence of some big theorems and one should prove the big theorems. I give only the bare outlines and references.

 The first step is to establish  a description of rational maps using linear systems. This is treated in many introductory algebraic geometry books, for example \cite[Chap.III]{shaf} or \cite[Chap.6]{MR0453732}. 
The  embedding  $X\into \p^{n+1}$ is given by the  linear system
of hyperplane sections, call it $|H|$. 

The hard theorem, called the Noether-Lefschetz theorem, says that if
$\dim X\geq 3$ then every linear system is a subsystem of some $|mH|$.
In Lefschetz's language this follows from
$H^1\bigl(X(\c), \z\bigr)=0$ and $H^2\bigl(X(\c), \z\bigr)\cong\z$;
 see \cite{MR0033557, SGA2} or \cite[p.156]{gri-har} for proofs. 
If a subsystem of  $|mH|$ gives an embedding then the degree of its image is $m^{\dim X}\cdot \deg X$. Thus  $X\into \p^{n+1}$ is the unique smallest degree embedding of $X$ into any projective space.
We are done if
$\dim X\geq 3$.
Another approach, which works also for surfaces 
will be discussed in Paragraph~\ref{isom.hyp.thm.pf.1}.

Finally, a theorem of Noether \cite[Sec.5]{MR1579926} 
describes all linear systems on plane curves that have unusually large dimension. A particular case of this 
says that an irreducible plane curve of degree $\geq 4$ has a unique embedding into $\p^2$. A complete proof is given by
Hartshorne  \cite{MR857224}, see also  \cite[p.56]{ACGH}.  
\end{say}

\section{Non-rationality of large degree hypersurfaces}\label{sec.7}

We start with the following answer to
Question~\ref{ques.2.1.s}.

\begin{thm}\label{lar.deg.nonrat.thm}
 Let $X\subset \p^{n+1}$ be a smooth hypersurface of degree $d$.
If $d\geq n+2$ then $X$ is not rational. More generally, the image of every 
rational map $\Phi:\p^n\map X$  has dimension $\leq n-1$.
\end{thm}

This will be a direct consequence of Propositions~\ref{non.unirat.prop}~and~\ref{vol.f.proj.say}, and
the same method also answers Question~\ref{ques.2.2.s}
for large degree hypersurfaces.

\begin{thm} \label{lar.deg.bir.thm}
Let $X_1, X_2\subset \p^{n+1}$ be  smooth hypersurfaces of degrees $d_1, d_2$ that are birational to each other. Assume that  $d_1\geq n+3$.

Then $d_1=d_2$ and the $X_i$ can be obtained from each other  by a linear change of coordinates.
\end{thm}

In order to prove Theorem~\ref{lar.deg.nonrat.thm}, we need to study algebraic differential forms.

\begin{defn} \label{alg.deff.forms}
Let $X$ be a smooth hypersurface of dimension $n$. 
As we noted in Definition~\ref{smooth.defn}, it is also a real manifold of dimension $2n$.
An {\it algebraic differential $m$-form}  is a differential $m$-form that locally can be written as linear combination of terms 
$$
\phi({\mathbf x})\cdot d\psi_1({\mathbf x})\wedge \cdots \wedge  d\psi_m({\mathbf x}),
\eqno{(\ref{alg.deff.forms}.1)}
$$
where $\phi$ and the $\psi_i$ are rational functions. An algebraic differential form is {\it defined} at a point $p\in X$ if it can be written with summands as in 
(\ref{alg.deff.forms}.1) where $\phi$ and the $\psi_i$ are all defined at $p$.

We will be especially interested in $n$-forms, which are
the  {\it algebraic volume forms.}  
(It is not hard to check that if $\sigma$  is an  algebraic volume form
then $\sqrt{-1}^{n}\sigma\wedge\bar{\sigma}$ is a usual volume form, albeit possibly with zeros and poles.)
\end{defn}

\begin{exmp}[Algebraic volume forms on $\c^n$]\label{vol.f.Cn.aff.say}
On $\c^n$ the ``standard'' volume form is
$dz_1\wedge\cdots\wedge dz_n$. If $\psi$ is a rational function then
$$
d\psi=\tsum_i \tfrac{\partial \psi}{\partial z_i}\ dz_i,
$$
hence  every algebraic volume form can be written as
$$
\sigma:=\phi({\mathbf z})\cdot dz_1\wedge\cdots\wedge dz_n,
\eqno{(\ref{vol.f.Cn.aff.say}.1)}
$$
where $\phi$ is a rational function. This in undefined precisely where $\phi$ is undefined. Thus, as we noted in Definition~\ref{rat.map.aff.defn}, either $\sigma$ is everywhere defined or it is undefined exactly along a hypersurface. 
\end{exmp}

The projective case behaves quite differently.

\begin{lem}\label{vol.f.Pn.aff.say}
There are no  everywhere defined algebraic volume forms on $\p^n$.
\end{lem}

Proof. 
We use projective coordinates $X_0{:}\cdots{:}X_n$ and 
show the claim using just 2 of the 
affine charts. Let these be 
$$
\begin{array}{lcl}
(x_1,\dots, x_n)&=&\bigl({X_1}/{X_0},\dots,   {X_n}/{X_0}\bigr)\\
(y_0,\dots, y_{n-1})&=&\bigl({X_0}/{X_n},\dots,  {X_{n-1}}/{X_n}\bigr).
\end{array}
$$
Thus we can transition between these charts by the formulas
$$
x_1=\tfrac{y_1}{y_0}, \dots, x_{n-1}=\tfrac{y_{n-1}}{y_0}, x_n=\tfrac1{y_0}.
$$
This gives that 
$$
dx_1\wedge\cdots\wedge dx_n=
d\tfrac{y_1}{y_0}\wedge \cdots\wedge d\tfrac{y_{n-1}}{y_0}\wedge d\tfrac1{y_0}=
\tfrac{(-1)^n}{y_0^{n+1}}\cdot dy_0\wedge\cdots\wedge dy_{n-1}.
$$
By Example~\ref{vol.f.Cn.aff.say}, 
any algebraic volume form on the first chart $\c^n_{\mathbf x}$ is of the form
$f({\mathbf x})\cdot dx_1\wedge\cdots\wedge dx_n$. On the second chart  $\c^n_{\mathbf y}$ it becomes
$$
\tfrac{(-1)^n}{y_0^{n+1}}
f\bigl(\tfrac{y_1}{y_0}, \dots, \tfrac{y_{n-1}}{y_0},\tfrac1{y_0}\bigr)
\cdot dy_0\wedge\cdots\wedge dy_{n-1}.
$$
Thus it has a pole of order $\geq n+1$ along $(y_0=0)$, 
so not defined there. \qed

\begin{prop}\label{non.unirat.prop}
  Let $X$ be a smooth projective hypersurface  that has an everywhere defined,  nonzero,  algebraic volume form $\sigma_X$. Then $X$ is not rational. In fact, there is not even a
rational map $\Phi:\p^n\map X$ with dense image. 
\end{prop}

Proof. By (\ref{rat.map.pr.defn}.6)  (or the more general Theorem~\ref{rtl.cod.2.thm}),  there is a closed subset $Z\subset \p^n$ of codimension $\geq 2$  such that $\Phi$ is defined on $\p^n\setminus Z$. Thus $\Phi^*\sigma_X$ is an algebraic volume form on
$\p^n\setminus Z$. By the purity of poles noted in Example~\ref{vol.f.Cn.aff.say}, 
$\Phi^*\sigma$ must extend to  an everywhere defined algebraic volume form on
$\p^n$. This is impossible by Lemma~\ref{vol.f.Pn.aff.say}. \qed
\medskip

({\it Side remark.} This is one result that is quite different in
positive characteristic. If we are in characteristic $p$ then
$(x^p)'=px^{p-1}$ is identically 0. Thus the pull-back of a
nonzero  algebraic volume form can be identically zero.)
\medskip

We should thus write down algebraic volume forms on  hypersurfaces. We start with the affine case.

\begin{say}[Algebraic volume forms on affine hypersurfaces]\label{vol.f.hyp.aff.say}
Let  $X\subset \c^{n+1}$ be a  hypersurface. As in \cite[Sec.III.6.4]{shaf}, it is quite easy to write down a nowhere zero volume form on $X^{\rm sm}$, the smooth part of $X$.
\medskip

{\it Claim \ref{vol.f.hyp.aff.say}.1.} For any hypersurface $X=(h=0)\subset \c^{n+1}$,  the volume forms
$$
\sigma_i:=\tfrac{(-1)^i}{\partial h/\partial z_i} dz_1\wedge\cdots\wedge\widehat{dz_i}\wedge\cdots\wedge dz_{n+1}
\eqno{(\ref{vol.f.hyp.aff.say}.2)}
$$
patch together to  a volume form $\sigma_X$ that is
defined and nowhere zero  on $X^{\rm sm}$. 
\medskip

Proof. By definition $h|_X=0$, thus 
$dh=\tsum_j \tfrac{\partial h}{\partial z_i}dz_j$ 
is also zero on $X$. Wedging it with all the $dz_i$ except
$dz_{i_1}$ and $dz_{i_2}$ gives that $\sigma_{i_1}=\sigma_{i_2}$. Thus
the $\sigma_i$ patch together to a rational volume form on $X$. 

It is clear that $\sigma_i$ is defined  and nowhere zero on the open set where    $\partial h/\partial z_i\neq 0$.
Thus at least one of the $\sigma_i$ is defined at a smooth point of $X$, so
$\sigma_X$ is defined  and nowhere zero  on  $X^{\rm sm}$. \qed

\medskip
As a counterpart/consequence of (\ref{rat.map.hyp.defn}.2) we obtain the following. 
\medskip

{\it Claim \ref{vol.f.hyp.aff.say}.3.} Let $X\subset \a^{n+1}$ be a smooth hypersurface and  $\rho$  an everywhere defined algebraic volume form on $X$. Then it can be written as $\rho=p({\mathbf x})\sigma_X$ where $p$ is a polynomial on $\a^{n+1}$. \qed

\medskip

We already noted the following analog of (\ref{rat.map.hyp.defn}.1) for $X=\c^n$, the general case is proved similarly using (\ref{rat.map.hyp.defn}.1).

\medskip

{\it Purity Principle \ref{vol.f.hyp.aff.say}.4.}  Let $\rho$ be an algebraic volume form on a smooth hypersurface $X$. Then either $\rho$ is everywhere defined or it is not defined along a codimension 1 subset of $X$.
\end{say}

\begin{prop} \label{vol.f.proj.say}
Let  $X\subset \p^{n+1}$ be a smooth hypersurface of degree $d$. If $d\geq n+2$ then it has a nonzero, everywhere defined algebraic volume form. 

More precisely,  
 the vector space of  everywhere defined algebraic volume forms  on $X$ is naturally isomorphic to the 
vector space of homogeneous polynomials of degree $d-n-2$ in
$n+2$ variables. 
\end{prop}

Proof. We can combine the computations of (\ref{vol.f.hyp.aff.say}) and (\ref{vol.f.Pn.aff.say}) to determine the volume forms on
smooth hypersurfaces, but the following way may be quicker.

As we discussed in Definition~\ref{proj.sp.defn},  $\c\p^{n+1}$ is the quotient of  $\c^{n+2}\setminus\{0\}$ by  the
$\c^*$-action
$$
m_{\lambda}:(X_0,\dots, X_{n+1})\mapsto (\lambda X_0,\dots, \lambda X_{n+1}),
\eqno{(\ref{vol.f.proj.say}.1)}
$$
which is induced by the vector field
$v_X:=\tsum X_i\tfrac{\partial}{\partial X_i}$.

We can thus think of a volume form on $\c\p^{n+1}$ as a 
$\c^*$-equivariant volume form on $\c^{n+2}\setminus\{0\}$,
contracted by the vector field  $v_X$.
By (\ref{vol.f.hyp.aff.say}.3) every volume form on $\c^{n+2}\setminus\{0\}$
extends to a volume form on $\c^{n+2}$.

Let next $X=X(H)\subset \c\p^{n+1}$ be a smooth hypersurface of degree $d$ and
$C_X:=(H=0)\subset \c^{n+2}$ the cone (\ref{proj.sp.defn}.5)
over it. Denote the cone minus its vertex by 
$C_X^{\circ}$. We can thus view a 
volume form on $X$  as a $\c^*$-equivariant volume form on $C_X^{\circ}$ contracted by $v_X$. 
We wrote down a nowhere zero volume form on $C_X^{\circ}$ in (\ref{vol.f.hyp.aff.say}.2). Choose the chart where it is given by 
$$
\Sigma_X:=\Bigl\{\tfrac{\pm 1}{\partial H/\partial X_{n+1}} dX_0\wedge\cdots\wedge dX_n\Bigr\}.
\eqno{(\ref{vol.f.proj.say}.2)}
$$
In particular,
$$
m_{\lambda}^*\Sigma_X=\lambda^{n+1-(d-1)}\Sigma_X. 
\eqno{(\ref{vol.f.proj.say}.3)}
$$
Thus if $G$ is any homogeneous polynomial then
$$
m_{\lambda}^* \bigl(G\Sigma_X)= \lambda^{\deg G+n+2-d}\cdot G\Sigma_X.
\eqno{(\ref{vol.f.proj.say}.4)}
$$
Hence $G\Sigma_X$  is $\c^*$-invariant iff $\deg G=d-n-2$.

Let us work out explicitly  what we get in therms of the  forms
$\sigma_i$ given in (\ref{vol.f.hyp.aff.say}.2).
Contracting $G\Sigma_X$ by $v_X$ gives
$$
\tsum X_i\tfrac{\partial}{\partial X_i}\bigl(G\Sigma_X)=
\tsum_{i=0}^n X_iG\cdot 
\frac{dX_0\wedge\cdots\wedge\widehat{dX_i}\wedge\cdots\wedge dX_n}{\partial H/\partial X_{n+1}}.
\eqno{(\ref{vol.f.proj.say}.5)}
$$
We need to pull this back to  $X$ by the map
$(x_1,\dots, x_{n+1})\mapsto  (1{:}x_1{:}\cdots{:} x_{n+1})$. 
The pull-back of $dX_0$ is then 0 and the only term that survives is
$$
 G(1, x_1,\dots, x_{n+1})
\tfrac{dx_1\wedge\cdots\wedge dx_n}
{\partial H(1, x_1,\dots, x_{n+1})/\partial x_{n+1}}=G(1, x_1,\dots, x_{n+1})\cdot \sigma_{n+1}.
\eqno{(\ref{vol.f.proj.say}.6)}
$$
where $\sigma_{n+1}$ is as  in (\ref{vol.f.hyp.aff.say}.2). \qed

\begin{say}[Comments on the proofs of Theorems~\ref{isom.hyp.thm} and \ref{lar.deg.bir.thm}]\label{isom.hyp.thm.pf}
We proved in Proposition~\ref{vol.f.proj.say} that the 
vector space of homogeneous polynomials of degree $\deg X-n-2$ in
$n+2$ variables  (an invariant of the embedding $X\into \p^{n+1}$) is naturally isomorphic to
 the vector space of  everywhere defined algebraic volume forms  on $X$ (an intrinsic invariant of $X$).

We are completely done with both theorems if $d=n+3$. Then the volume forms tell us the linear polynomials and they define the embedding $X\into \p^{n+1}$. 

We are almost done if  $d>n+3$; the  missing ingredient is that
$H^2\bigl(X(\c), \z\bigr)$ is torsion free for $\dim X\geq 2$. This is another special case  of the Noether-Lefschetz theorem.

If $d<n+3$ then we can use the duals of holomorphic volume forms.
These are wedge products of tangent vector fields. We can pull back a vector field by an isomorphism but not by an arbitrary map.
This is why Theorem~\ref{isom.hyp.thm}  works for $d<n+3$ but 
 Theorem~\ref{lar.deg.bir.thm} does not.

For (\ref{lar.deg.bir.thm}.1) we need to combine these with the methods of
Paragraph~\ref{isom.hyp.thm.pf.1}.
\qed
\end{say}

\section{Non-rationality of low degree hypersurfaces}\label{sec.8}

In the previous Section we  proved that  a smooth hypersurface
$X\subset \p^{n+1}$ of degree $\geq n+2$ is not rational.
Here we discuss what is know when the degree is $\leq n+1$. 
Most of the results use different---and much more advanced---methods, so
we  give the  statements only, with barely a hint of how they can be proved. 

We proved  in Theorem~\ref{stere.gen.quad.thm}  that every irreducible quadric hypersurface is rational over $\c$ and, by Corollary~\ref{cubic.C.rtl.prop},  
 every smooth cubic surface is rational.

We already noted that    smooth cubic 3-folds are all   non-rational
 and the same holds for smooth quartic 3-folds \cite{MR0291172};
we say more about this in  Section~\ref{sec.10}.

Another method, using differential forms over fields of positive characteristic was introduced in
\cite{Kollar95a}. We proved in (\ref{vol.f.proj.say}.5) that if $\deg X<n+2$ then there are no 
everywhere defined $n$-forms on $X$. One can show that there are also no
everywhere defined algebraic $r$-forms on $X$ for $r\leq n$.

It turns out that while the first of these claims continues to hold over fields of positive characteristic, there are some hypersurfaces that carry $(n-1)$-forms. There are some  technical issues with singularities of these hypersurfaces, but one can sometimes conclude that these  hypersurfaces are non-rational. However, this happens in positive characteristic. Nonetheless, one can use some general theorems going back to Matsusaka and Mumford \cite{mats-mumf}  to get similar conclusions over $\c$.

A serious drawback is that we can prove non-rationality only for hypersurfaces
$X=(\sum_Ia_Ix^I=0)$ whose coefficients satisfy countably many  
conditions of the form  $p_j(a_I)\neq 0$ where the $p_j$ are (pretty much unknown) polynomials. 
We refer to  such hypersurfaces as {\it very general.}

\begin{thm} \cite{Kollar95a}  Very general hypersurfaces $X_d\subset \p^{n+1}$  of degree $d\geq \tfrac23 n+3$ are not rational. 
Moreover, $X_d$ is not birational to
any product  $\p^1\times Y$ where $Y$ is a hypersurface (or variety) of dimension $n-1$. \end{thm}

I stress that we really do not know which hypersurfaces are covered by the theorem, though some concrete examples were written down in \cite[Sec.V.5]{rc-book}. 

Another argument, this time using degeneration to certain singular
varieties and the topology of their resolution was introduced by
Voisin \cite{MR3359052}. It was used by many authors to prove non-rationality of  varieties. For hypersurfaces the strongest result is proved by \cite{schreieder}.

\begin{thm}[Schreieder]  Very general hypersurfaces $X_d\subset \p^{n+1}$  of degree $d\geq 3+\log_2n$ are not rational.
\end{thm}

\section{Rigidity of low degree  hypersurfaces}\label{sec.10}

Much less is known about Question~\ref{ques.2.2} for
hypersurfaces $X_d\subset \p^{n+1}$ whose degree is $3\leq d\leq n+1$.
We studied the degree 3 case in Section~\ref{sec.5}. 
The answer in the $d=n+1$ case is 
one of the crowning achievements of  the Noether--Fano rigidity theory.

I call it the {N}oether--{F}ano--{S}egre--{I}skovskikh--{M}anin--{C}orti--{P}ukhlikov--{C}heltsov--de{\,}{F}ernex--{E}in--{M}usta{\cedilla{t}\u{a}}--{Z}huang theorem,
although maybe {N}oether--{F}ano--{S}egre--{I}skovskikh-{M}anin--{I}skovskikh--{C}orti--{P}ukhlikov--{C}heltsov--{P}ukhlikov--de{\,}{F}ernex-{E}in-{M}usta{\cedilla{t}\u{a}}--de{\,}{F}ernex--{Z}huang theorem
would be a historically more accurate name to emphasize the especially large contributions of {I}skovskikh and {P}ukhlikov.  See \cite{k-ah2} for a detailed survey.

\begin{thm}\label{superrigid.thm}
 Let $X\subset   \p^{n+1}$ be a smooth hypersurface of degree $n+1$. Assume that $n\geq 3$.  Then $X$ is not birational to any other smooth hypersurface.
\end{thm}

\begin{say}[A short history of Theorem~\ref{superrigid.thm}] 

The first similar result is Max~Noether's description of all birational maps $\p^2\map \p^2$ \cite{MR1509694}, whose method formed the basis of all further developments. 

 Theorem~\ref{superrigid.thm} was first stated by
Fano for 3-folds \cite{Fano-1908, Fano-1915}. His arguments contain many of the key ideas, but they also have  gaps. 
I call this approach the {\it   Noether--Fano method.}
The first complete proof  for 3-folds, along the lines indicated by Fano, is in
 Iskovskikh-Manin \cite{MR0291172}.  Iskovskikh and his school used this method to prove similar results for 
many other 3-folds, see \cite{MR537686,  Sar81, MR1668579, MR1859707}.
This approach  was gradually  extended to higher dimensions  by
Pukhlikov  \cite{MR870730, MR1650332, MR1970356} and  Cheltsov \cite{MR1805602}.
These results were complete up to dimension 8, but needed some additional general position assumptions in higher dimensions.
 A detailed  survey of this direction  is  in  \cite{puk-book}.

The  Noether--Fano method and the  Minimal Model Program  were brought  together   by 
Corti \cite{MR1311348}. 
(See \cite{k-oldintro} for a by now outdated but elementary introduction and \cite{km-book} for a detailed treatment.)
Corti's technique  has been very successful in many cases, especially for 3--folds; see \cite{MR1798978} for a detailed study and \cite[Chap.5]{ksc} for an introduction. 
However, usually one needs some special tricks to make the last steps work, and a good higher dimensional version proved elusive for a long time.

New methods involving multiplier ideals were introduced by 
 de~Fernex-Ein-Musta\cedilla{t}\u{a} \cite{MR1981899}; these led to a more streamlined proof that worked up to dimension 12.   
The proof of Theorem~\ref{superrigid.thm} was finally completed by de~Fernex \cite{MR3455160}.

The recent paper of Zhuang \cite{zhuang}  makes the final step of  the Corti approach much easier in higher dimensions.
 The papers \cite{sti-zhu, zhuang, liu-zhuang} contain more general results and  applications. 
\end{say}

 \section{Connections with  the classification of varieties}

So far I have studiously avoided assuming prior knowledge of algebraic geometry.
However, in order to explain the place of these results in algebraic geometry,  it becomes necessary to use the basic theory, as in  \cite[Vol.I]{shaf}. 
In particular we need to  know the notions of  ampleness and  canonical class  $K_X$.

The classification theory of algebraic varieties---developed by Enriques for surfaces and extended by Iitaka and then Mori to higher dimensions---says that every variety can be built from 3 basic types:
\begin{itemize}
\item (General type)  $K_X$ is ample,
\item (Calabi-Yau)   $K_X$ is trivial and
\item (Fano)  $-K_X$ is ample;
\end{itemize}
see \cite{k-icm} for an introduction. 
Moreover, in the Fano case the truly basic ones are those that have
{\it class number} equal to 1. That is, every divisor $D$ on $X$ is
linearly equivalent to a (possibly rational) multiple of $-K_X$. 
(For smooth varieties, the class number is the same as the {\it Picard number.})
Such examples are $\p^n$  (where every divisor is linearly equivalent to 
$\frac{m}{n+1}(-K_{\p^n})$ for some $m\in \z$) or smooth hypersurfaces 
$X\subset \p^{n+1}$ of
dimension $n\geq 3$ and of degree $d\leq n+1$ (where every divisor is linearly equivalent to 
$\frac{m}{n+2-d}(-K_X)$ for some $m\in \z$).

The computations of Proposition~\ref{vol.f.proj.say} show that if $X\subset \p^{n+1}$ is a smooth hypersurface then it is of general type iff $\deg X\geq n+3$,  Calabi-Yau   iff $\deg X= n+2$ and Fano   iff $\deg X\leq n+1$.

As in Theorem~\ref{lar.deg.bir.thm}, one can see  that if 2 varieties $X_1, X_2$  on the basic type list  list are birationally equivalent then they have the same type. Furthermore,
in the general type case they are even isomorphic; this generalizes Theorem~\ref{lar.deg.bir.thm}. 
In the Calabi-Yau case  $X_1$ and  $X_2$ need not be   isomorphic,  but the 
possible  birational maps between two  Calabi-Yau varieties  are reasonably well understood, especially in dimension 3;  see, for example,  \cite{MR986434, k-etc, MR3504536}. 

It is very useful to think of Fano varieties as representatives of 
{\it rationally connected varieties,} but this would have taken us in another direction; see
\cite{MR1848255} for an overview and \cite{rc-book, ar-ko} for more detailed treatments.

Thus, from the general point of view, for low degree cases, the best variant of 
Question~\ref{ques.2.2} is the following.

\begin{defn}  A Fano variety $X$ with class number 1 is called {\it weakly superrigid} if 
every birational map  $\Phi:X\map Y$ to another Fano variety $Y$ with class number 1  is an isomorphism.
\end{defn}

\begin{ques}\label{ques.2.2.v} Which Fano varieties are weakly superrigid?
\end{ques}

 The adjective ``weakly'' is not standard, I use it just to avoid further definitions. 
The correct definition of {\it  superrigid} allows  $Y$ to have terminal singularities and to be  a Mori fiber space; see \cite{MR1730270, MR2195677}.  
The proofs in the Noether--Fano  theory are designed to prove superrigidity.
However, there should be many varieties that are weakly superrigid but not
superrigid, see Question~\ref{weak.neq.sup.hyp.ques}.

There are many Fano varieties $X$, especially in dimensions 2 and 3,
that are not birational to any other Fano variety but there are
 birational maps  $\Phi:X\map X$ that are not isomorphisms.
Such Fano varieties are called {\it rigid.}
It seems to me that superrigidity is the more basic notion, though, in dimension 3, the theory of rigid Fano varieties is very rich.

\section{Open problems about hypersurfaces}

The following  questions are stated in the strongest forms that are consistent with the known examples. I have no reasons to believe that the answer to any of them is positive and there may well be rather simple counter examples. As far as I know, there has been very little work on low degree hypersurfaces beyond cubics in dimension 4.

All the known rationality constructions work for even dimensional cubics only, so the following is still open.

\begin{ques} Is there any odd dimensional, smooth, rational cubic hypersurface?
\end{ques}

We leave it  to the reader the easy task of  generalizing Proposition~\ref{cubic.2L.rtl.prop}
 to the following higher dimensional version. 

\begin{prop}\label{cubic.2n.2L.rtl.prop}
 Let $X\subset \p^{2n+1}$ be a cubic hypersurface that contains 2 disjoint linear subspaces of dimension $n$. Then $X$ is rational. \qed
\end{prop}

This leaves the following open.

\begin{ques} Is the general even  dimensional, smooth,  cubic hypersurface
of dimension $\geq 4$  non-rational?
\end{ques}

 We saw that quadrics are rational and
in Section \ref{sec.5} we gave several examples of cubic hypersurfaces that are rational. However, no rationality construction is known for
smooth hypersurfaces of higher degree. 

 Example~\ref{rtl.sing.exmps}.2 gives higher degree rational
hypersurfaces $X_{2d+1}^{2n}\subset \p^{2n+1}$, but they are always singular. 
However,  the singularities are mild
when $d\leq n$. To be precise, they are   canonical if $d\leq n$ and terminal if $d\leq n-1$.  (See
\cite[Sec.2.3]{km-book} or \cite{kk-singbook} for  introductions to such singularities.)
Thus the following problem is still open.  

\begin{ques} Is there any  smooth, rational  hypersurface of degree $\geq 4$? 
\end{ques}

Understanding birational maps between Fano hypersurfaces
is even harder.  Theorem~\ref{superrigid.thm} deals with smooth hypersurfaces $X\subset \p^{n+1}$ of degree $n+1$.

\begin{ques} \label{weak.neq.sup.hyp.ques} Is every birational map between
 smooth hypersurfaces of degree $\geq 5$  an isomorphism?
\end{ques}

Here $\geq 5$ is necessary since 
there are  some smooth quartics with nontrivial birational maps.

\begin{exmp}\label{quart.2n.2L.inv.prop}
 Let $X\subset \p^{2n+1}$ be a quartic hypersurface that contains 2 disjoint linear subspaces $L_1, L_2$ of dimension $n$. 

As in Proposition~\ref{cubic.2L.rtl.prop}, for every ${\mathbf p}\in \p^{2n+1}\setminus (L_1\cup L_2)$ there is a unique line $\ell_{\mathbf p}$ through ${\mathbf p}$ that meets both $L_1, L_2$. 
This line meets $X$ in 4 points, two of these are on $L_1, L_2$. If ${\mathbf p}\in X$ then this leaves a unique 4th intersection point, call it $\Phi({\mathbf p})$. 
Clearly $\Phi$ is an involution which is not defined at ${\mathbf p}$ if either
${\mathbf p}\in L_1\cup L_2$ or if $\ell_{\mathbf p}\subset X$. 

To get a concrete example, 
for $n\geq 2$ consider  the smooth, quartic hypersurface
$$
X:=\bigl(x_0^4+\cdots + x_{2n+1}^4=0\bigr)\subset \p^{2n+1}.
$$
Then $\aut(X)$ is finite (probably of order  $4^{2n+1}\cdot (2n+2)!$) but combining 
Example~\ref{lines.on.diag.exmp} with the above observation 
 gives many birational involutions on $X$.
Most likely these involutions generate an infinite subgroup of $\bir(X)$. 
\end{exmp}


\def\cprime{$'$} \def\cprime{$'$} \def\cprime{$'$} \def\cprime{$'$}
  \def\cprime{$'$} \def\cprime{$'$} \def\dbar{\leavevmode\hbox to
  0pt{\hskip.2ex \accent"16\hss}d} \def\cprime{$'$} \def\cprime{$'$}
  \def\polhk#1{\setbox0=\hbox{#1}{\ooalign{\hidewidth
  \lower1.5ex\hbox{`}\hidewidth\crcr\unhbox0}}} \def\cprime{$'$}
  \def\cprime{$'$} \def\cprime{$'$} \def\cprime{$'$}
  \def\polhk#1{\setbox0=\hbox{#1}{\ooalign{\hidewidth
  \lower1.5ex\hbox{`}\hidewidth\crcr\unhbox0}}} \def\cdprime{$''$}
  \def\cprime{$'$} \def\cprime{$'$} \def\cprime{$'$} \def\cprime{$'$}
\providecommand{\bysame}{\leavevmode\hbox to3em{\hrulefill}\thinspace}
\providecommand{\MR}{\relax\ifhmode\unskip\space\fi MR }
\providecommand{\MRhref}[2]{%
  \href{http://www.ams.org/mathscinet-getitem?mr=#1}{#2}
}
\providecommand{\href}[2]{#2}

\bigskip

\noindent  Princeton University, Princeton NJ 08544-1000

{\begin{verbatim} kollar@math.princeton.edu\end{verbatim}}

\end{document}